\documentclass[journal]{IEEEtran}

\usepackage{algorithmic}
\usepackage[cmex10]{amsmath}
\usepackage{amssymb}
\usepackage{array}
\usepackage{bm}
\usepackage{cite}
\usepackage{float}
\usepackage{mathtools}
\usepackage{pgfplots}
\usepackage{pgfplotstable}
\usepackage{tikzscale}
\usepackage{url}

\usepgfplotslibrary{colorbrewer}
\pgfplotsset{cycle list/Set1-9}
\pgfplotsset{compat=1.16}

\ifCLASSOPTIONcompsoc
  \usepackage[caption=false,font=normalsize,labelfont=sf,textfont=sf]{subfig}
\else
  \usepackage[caption=false,font=footnotesize]{subfig}
\fi

\makeatletter
\let\old@ps@headings\ps@headings
\let\old@ps@IEEEtitlepagestyle\ps@IEEEtitlepagestyle
\def\psccfooter#1{%
    \def\ps@headings{%
        \old@ps@headings%
        \def\@oddfoot{\strut\hfill#1\hfill\strut}%
        \def\@evenfoot{\strut\hfill#1\hfill\strut}%
    }%
    \def\ps@IEEEtitlepagestyle{%
        \old@ps@IEEEtitlepagestyle%
        \def\@oddfoot{\strut\hfill#1\hfill\strut}%
        \def\@evenfoot{\strut\hfill#1\hfill\strut}%
    }%
    \ps@headings%
}
\makeatother

\floatstyle{ruled}
\newfloat{model}{thp}{lop}
\floatname{model}{Model}

\begin{document}

\title{Convex Relaxations of Maximal Load Delivery for Multi-contingency Analysis of Joint Electric Power and Natural Gas Transmission Networks}

\author{Byron~Tasseff,
        Carleton~Coffrin,
        and~Russell~Bent%
\thanks{B. Tasseff is with the Department
of Industrial and Operations Engineering, University of Michigan, Ann Arbor,
MI, 48109 USA, and Los Alamos National Laboratory, Los Alamos, NM, 87545 USA, e-mail: btasseff@lanl.gov.}%
\thanks{C. Coffrin and R. Bent are also with Los Alamos National Laboratory.}}%

\maketitle

\begin{abstract}
Recent increases in gas-fired power generation have engendered increased interdependencies between natural gas and power transmission systems.
These interdependencies have amplified existing vulnerabilities to gas and power grids, where disruptions can require the curtailment of load in one or both systems.
Although typically operated independently, coordination of these systems during severe disruptions can allow for targeted delivery to lifeline services, including gas delivery for residential heating and power delivery for critical facilities.
To address the challenge of estimating maximum joint network capacities under such disruptions, we consider the task of determining feasible steady-state operating points for severely damaged systems while ensuring the maximal delivery of gas and power loads simultaneously, represented mathematically as the nonconvex joint Maximal Load Delivery (MLD) problem.
To increase its tractability, we present a mixed-integer convex relaxation of the MLD problem.
Then, to demonstrate the relaxation's effectiveness in determining bounds on network capacities, exact and relaxed MLD formulations are compared across various multi-contingency scenarios on nine joint networks ranging in size from $25$ to $1{,}191$ nodes.
The relaxation-based methodology is observed to accurately and efficiently estimate the impacts of severe joint network disruptions, often converging to the relaxed MLD problem's globally optimal solution within ten seconds.
\end{abstract}


\IEEEpeerreviewmaketitle

\section*{Nomenclature}
\begin{IEEEdescription}[\IEEEusemathlabelsep\IEEEsetlabelwidth{$Y^s = g^s-$}]
    \item [\textbf{Electric Power Sets}]
    \item [{$\mathcal{N}$}] Set of buses
    \item [{$\mathcal{E}$}] Set of forward-oriented lines (branches)
    \item [{$\mathcal{E}^{R}$}] Set of reverse-oriented lines (branches)
    \item [{$\mathcal{G}$}] Set of bus-connected generators
    \item [{$\mathcal{L}$}] Set of bus-connected loads
    \item [{$\mathcal{H}$}] Set of bus-connected shunts

    \item [\textbf{Electric Power Parameters}]
    \item [{$\underline{V}_{i}, \overline{V}_{i} \geq 0$}] Lower, upper voltage mag. bounds for $i \in \mathcal{N}$
    \item [{$Y_{ij}, Y_{ij}^{c} \in \mathbb{C}$}] Line admittance, charging for $(i, j) \in \mathcal{E} \cup \mathcal{E}^{R}$
    \item [{$T_{ij} \in \mathbb{C}$}] Transformer properties along $(i, j) \in \mathcal{E}$
    \item [{$S_{i}^{d} \in \mathbb{C}$}] Maximum deliverable power at $i \in \mathcal{L}$
    \item [{$\underline{S}_{i}^{g}, \overline{S}_{i}^{g} \in \mathbb{C}$}] Lower, upper generation bounds for $i \in \mathcal{G}$
    \item [{$Y_{i}^{s}$}] Admittance of bus shunt $i \in \mathcal{H}$
    \item [{$\overline{S}_{ij}$}] Apparent power limit for $(i, j) \in \mathcal{E} \cup \mathcal{E}^{R}$
    \item [{$\underline{\theta}_{ij}^{\Delta}, \overline{\theta}_{ij}^{\Delta}$}] Lower, upper phase angle diffs. for $(i, j) \in \mathcal{E}$ \\

    \item [\textbf{Electric Power Variables}]
    \item [{$V_{i} \in \mathbb{C}$}] Voltage at bus $i \in \mathcal{N}$
    \item [{$S_{i}^{g} \in \mathbb{C}$}] Power supplied by $i \in \mathcal{G}$
    \item [{$S_{ij} \in \mathbb{C}$}] Power across line $(i, j) \in \mathcal{E} \cup \mathcal{E}^{R}$
    \item [{$z_{i}^{d} \in [0, 1]$}] Scalar for deliverable load at $i \in \mathcal{L}$
    \item [{$z_{i}^{s} \in [0, 1]$}] Scalar for fixed bus shunt at $i \in \mathcal{H}$
    \item [{$z_{i}^{v} \in \{0, 1\}$}] Energization status of bus $i \in \mathcal{N}$
    \item [{$z_{i}^{g} \in \{0, 1\}$}] Dispatch status of generator $i \in \mathcal{G}$

    \item [\textbf{Natural Gas Sets}]
    \item [{$\mathcal{J}$}] Set of junctions
    \item [{$\mathcal{R}$}] Set of receipts (producers)
    \item [{$\mathcal{D}$}] Set of deliveries (consumers)
    \item [{$\mathcal{A}$}] Set of junction-connecting components
    \item [{$\mathcal{P} \subset \mathcal{A}$}] Set of horizontal pipes
    \item [{$\mathcal{S} \subset \mathcal{A}$}] Set of short pipes
    \item [{$\mathcal{T} \subset \mathcal{A}$}] Set of resistors
    \item [{$\mathcal{V} \subset \mathcal{A}$}] Set of valves
    \item [{$\mathcal{W} \subset \mathcal{A}$}] Set of pressure-reducing regulators
    \item [{$\mathcal{C} \subset \mathcal{A}$}] Set of compressors
    \item [{$\delta_{i}^{+} \subset \mathcal{A}$}] Components directed \emph{from} $i \in \mathcal{J}$
    \item [{$\delta_{i}^{-} \subset \mathcal{A}$}] Components directed \emph{to} $i \in \mathcal{J}$

    \item [\textbf{Natural Gas Parameters}]
    \item [{$\underline{f}_{ij}, \overline{f}_{ij} \in \mathbb{R}$}] Lower, upper mass flow bounds for $(i, j) \in \mathcal{A}$
    \item [{$\underline{p}_{i}, \overline{p}_{i} \geq 0$}] Lower, upper pressure bounds for $i \in \mathcal{J}$
    \item [{$\overline{s}_{i} \geq 0$}] Upper supply mass flow bound for $i \in \mathcal{R}$
    \item [{$\overline{d}_{i} \geq 0$}] Upper demand mass flow bound for $i \in \mathcal{D}$
    \item [{$w_{ij} \geq 0$}] Resistance of pipe $(i, j) \in \mathcal{P}$
    \item [{$\tau_{ij} \geq 0$}] Resistance of resistor $(i, j) \in \mathcal{T}$
    \item [{$\underline{\alpha}_{ij}, \overline{\alpha}_{ij} \geq 0$}] Lower, upper scalars for $(i, j) \in \mathcal{W} \cup \mathcal{C}$

    \item [\textbf{Natural Gas Variables}]
    \item [{$f_{ij} \in \mathbb{R}$}] Mass flow along $(i, j) \in \mathcal{A}$
    \item [{$y_{ij} \in \{0, 1\}$}] Flow direction along $(i, j) \in \mathcal{A}$
    \item [{$s_{i} \geq 0$}] Supply at receipt $i \in \mathcal{R}$
    \item [{$d_{i} \geq 0$}] Demand at delivery $i \in \mathcal{D}$
    \item [{$p_{i} \geq 0$}] Pressure at junction $i \in \mathcal{J}$
    \item [{$z_{ij} \in \{0, 1\}$}] Status of controllable element $(i, j) \in \mathcal{V} \cup \mathcal{W}$

    \item [\textbf{Interdependency Modeling}]
    \item [{$\mathcal{K}$}] Set of links between $i \in \mathcal{D}$ and $j \in \mathcal{G}$
    \item [{$h_{ij}^{1, 2, 3}$}] Coefficients of heat rate curve for $(i, j) \in \mathcal{K}$
\end{IEEEdescription}

\section{Introduction}
\label{section:introduction}
\IEEEPARstart{B}{etween} 2012 and 2040, global electric power generation capacity is predicted to increase from $21.6$ million gigawatt-hours (GWh) to $36.5$ million GWh.
Of this, gas-fired generation is expected to increase from $22\%$ to $28\%$ \cite{conti2016international}.
This growing dependence underscores the increasing sensitivity of power systems to upstream disruptions in gas pipelines.
The most recent example is the February 2021 Texas power crisis, where the Electric Reliability Council of Texas experienced a loss of nearly $52.3$ GW ($48.6\%$) of its generation capacity.
Nearly half of the loss was attributed to a lack of gas-fired power generation \cite{ercot-slides}.
Other examples include the 2014 polar vortex, where curtailments in gas delivery resulted in roughly $25\%$ of generation outages throughout the Pennsylvania-New Jersey-Maryland Interconnection \cite{pjm2014}.
Disruptions in the gas grid can also inhibit the transport of fuel required for residential heating.
This begets an important tradeoff between gas delivery and power delivery during severe network disruptions.
Understanding these interdependencies is critical for the resilience of gas and power delivery systems.

\begin{figure}[t]
    \begin{center}
        \begin{tikzpicture}
    \begin{axis}[xlabel=Event,ylabel=Gas load delivered ($\%$),
                 domain=0:110,ymin=30,ymax=105,xmin=0,xmax=110,samples=1000,
                 xtick pos=bottom,xtick={15,40,55,90,100},xticklabels={(i),(ii),(iii),(iv),(v)},
                 width=0.85\linewidth,height=5.5cm,enlargelimits=false]
                 \draw[fill=gray!50,opacity=0.25,draw=none] plot[smooth, samples=100, domain=40:90]
                     (40,0) -| (40,105) -| (90,105) -| (90,0) -- cycle;
                 \draw[solid] plot [smooth, tension=0.0] coordinates {(0,100) (15,100)};
                 \draw[solid] plot [smooth, tension=1.0] coordinates {(15,100) (40,75) (70,100)};
                 \draw[solid] plot [smooth, tension=0.0] coordinates {(70,100) (110, 100)};
                 \draw[dotted] (15,34) -- (15,105);
                 \draw[dotted] (40,69) -- (40,105);
                 \draw[dotted] (55,69) -- (55,105);
                 \draw[dotted] (90,34) -- (90,105);
                 \draw[dotted] (100,34) -- (100,105);
    \end{axis}
    \begin{axis}[ylabel=Power load delivered ($\%$),axis y line*=right,
                 axis x line=none,
                 domain=0:110,ymin=30,ymax=105,xmin=0,xmax=110,samples=1000,
                 width=0.85\linewidth,height=5.5cm,enlargelimits=false]
                 \draw[dashed] plot [smooth, tension=0.0] coordinates {(0,100) (15,100)};
                 \draw[dashed] plot [smooth, tension=0.5] coordinates {(15,100) (20,98) (39,85) (60,90) (90,100)};
                 \draw[dashed] plot [smooth, tension=0.0] coordinates {(90,100) (110, 100)};
                 \node[align=left,anchor=west] (A) at (14.5,64.5) {\footnotesize (i) {Disruptive event begins.}};
                 \node[align=left,anchor=west] (B) at (14.5,55) {\footnotesize (ii) {Event and cascading effects}\\[-0.55em]
                                                               \footnotesize {cease. Load restoration processes}\\[-0.55em]
                                                               \footnotesize {begin without network repairs.}};
                 \node[align=left,anchor=west] (C) at (14.5,45.5) {\footnotesize (iii) {Network repairs commence.}};
                 \node[align=left,anchor=west] (D) at (14.5,40.4) {\footnotesize (iv) {Load restoration is complete.}};
                 \node[align=left,anchor=west] (E) at (14.5,35.5) {\footnotesize (v) {Network repairs are complete.}};
                 \node[align=left,anchor=east] (G) at (41,97){\footnotesize{Power}};
                 \node[align=left,anchor=east] (P) at (41,72) {\footnotesize{Gas}};
    \end{axis}
\end{tikzpicture}
    \end{center}
    \caption{A high-level illustration of natural gas and power transmission network responses to a hypothetical severe disruption. The shaded region indicates the points in the disruption and restoration timeline that are studied in this paper using an optimization-based assessment of damaged network capacities.}
    \label{figure:restoration}
\end{figure}
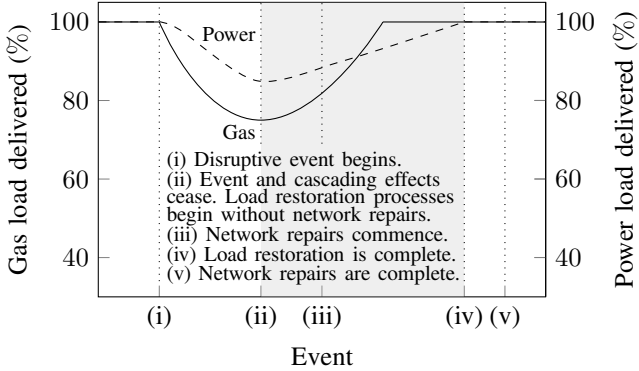

The contingency response measures considered in this paper are illustrated in Figure \ref{figure:restoration}.
Given a severe disruption, (i) gas and/or power load deliveries decrease as gas and/or power network elements are impaired and effects begin to cascade, and (ii) cascading effects subside, and a new stable operating point is realized.
After (ii), load can gradually be restored via operational methods until (iii) network repairs begin.
These restorative actions are performed until (iv) all gas and power loads can be delivered.
Repairs continue until (v) all gas and power network components are again operational.
Addressing all event types within Figure \ref{figure:restoration}, however, is a substantial task.
To make the scope more manageable, we focus more narrowly on ascertaining optimal steady-state operating points between events of types (ii) and (iv), i.e., decisions that maximize gas and power load delivery in the surviving gas-power system.

In this paper, this task is formalized as the steady-state joint Maximal Load Delivery (MLD) problem.
The problem is informally stated as follows: given severely damaged gas and power networks in which multiple components have become nonoperational, maximize the amounts of prioritized gas and active power loads that can be served simultaneously in the damaged joint network, subject to steady-state natural gas and alternating current (AC) power network physics.
The nonconvex physics and discrete nature of operations in the joint network (e.g., the opening and closing of valves in the gas network) render this a challenging mixed-integer nonlinear program (MINLP).
To increase its tractability, we develop a mixed-integer convex programming (MICP) relaxation of the MLD problem.
The MICP is found to be an effective means for bounding maximum total deliverable gas and power loads.

This paper expands upon existing MLD methods for independent gas \cite{tasseff2020natural} and power networks \cite{coffrin2018relaxations}, as well as approaches from joint network modeling, to formulate and solve the joint gas-power MLD problem.
Its contributions include
\begin{itemize}
    \item The first formulation of the gas-power MLD problem;
    \item A reliable MICP relaxation of the MLD problem;
    \item Proof-of-concept analyses of MLD gas-power tradeoffs.
\end{itemize}

The remainder of this paper proceeds as follows: Section \ref{section:background} reviews relevant gas, power, and joint steady-state optimization models that appear in the literature, then formulates the requirements for AC power and gas pipeline operational feasibility as an MINLP;
Section \ref{section:mld-formulations} formulates the MLD problem as an MINLP, then proposes an MICP relaxation;
Section \ref{section:computational-experiments} rigorously benchmarks the MINLP and MICP formulations across multiple joint gas-power networks of various sizes, then provides proofs of concept for joint multi-contingency analysis using the MLD method;
and Section \ref{section:conclusion} concludes the paper.

\section{Background for Network Modeling}
\label{section:background}
The past decade has seen remarkable theoretical and algorithmic advances in the independent fields of power and natural gas network optimization.
A recent survey of relaxations and approximations used in power system optimization is presented by \cite{molzahn2019survey}.
The study in power most related to this paper is by \cite{coffrin2018relaxations}, who introduce the AC MLD problem and propose various relaxations to increase its tractability.
The method was later extended by \cite{sundar2018probabilistic} to identify the $k$ components that maximize network disruption, as well as by \cite{sundar2019nk} to identify multi-contingency scenarios that would benefit from more detailed cascading analyses.
The MLD problem was also exploited by \cite{rhodes2020balancing}, who applied it within a bilevel optimization for balancing wildfire risk and power outages.
Finally, an implementation of the power MLD problem is presented by \cite{rhodes2021powermodelsrestoration}, who provide formulations via the \textsc{PowerModelsRestoration} package.
Their implementation, in fact, serves as a computational foundation for the power system modeling portion of this paper.

As with power, the growing utilization of gas networks has led to a variety of optimization studies.
A summary of recent work related to the optimization-based assessment of gas network capacities is provided by \cite{HILLER2018797}.
Steady-state models and approximations of gas network components amenable to optimization applications are provided by \cite{schmidt2016high}.
However, the study in gas most related to this paper is by \cite{tasseff2020natural}, who develop the steady-state gas MLD problem and an MICP relaxation.
Another related study is by \cite{ahumada2021nk}, who examine the the problem of identifying the $k$ components of a gas network whose simultaneous failure maximizes disruption to the network.

An even more recent body of literature has examined the optimal coordination of gas and power infrastructures.
A review of joint gas and power planning is given by \cite{farrokhifar2020energy}.
Other studies have focused on market coordination and energy pricing problems \cite{byeon2019unit,mitridati2020heat}.
Many studies have assumed the networks to be fully coordinated, examining optimal scheduling of generator dispatching and gas compressor operations \cite{zlotnik2016coordinated}.
Recent studies have expanded upon these earlier joint ``optimal gas-power flow'' problems, developing specialized formulations and algorithms for related applications \cite{wang2017convex,jiang2018coordinated,mirzaei2020novel}.
A smaller number of studies have considered joint problems related to restoration, e.g., scheduling of general large-scale interdependent infrastructures in \cite{abeliuk2016interdependent} and last-mile restoration of joint gas and power systems in \cite{coffrin2012last}.
The remaining subsections build upon previous studies to define the requirements for steady-state operation of a damaged joint gas-power network.

\subsection{Power Transmission Network Modeling}
\label{subsection:power-modeling}

\paragraph{Notation for Sets}
A power network is represented by an arbitrarily directed graph $(\mathcal{N}, \mathcal{E} \cup \mathcal{E}^{R})$, where $\mathcal{N}$ is the set of buses, $\mathcal{E}$ is the set of forward-directed branches (or lines), and $\mathcal{E}^{R}$ is the set of branches in their reverse orientation.
The set of generators (producers), loads (consumers), and shunts are denoted by $\mathcal{G}$, $\mathcal{L}$, and $\mathcal{H}$, respectively, which are attached to existing buses $i \in \mathcal{N}$.
We let the subset of these components attached to $i \in \mathcal{N}$ be denoted by $\mathcal{G}_{i}$, $\mathcal{L}_{i}$, and $\mathcal{H}_{i}$.
We next define the decision variables and constraints required to model a damaged AC power network's steady-state operations.

\paragraph{Power Network Modeling Requirements}
\begin{model}[t]
\begin{subequations}
    \begin{align}
        & S_{ij} = \left(Y_{ij} + Y^{c}_{ij} \right)^{*} \frac{\lvert V_{i} \rvert^{2}}{\lvert T_{ij} \rvert^{2}} - Y^{*}_{ij} \frac{V_{i} V_{j}^{*}}{T_{ij}}, ~ \forall (i, j) \in \mathcal{E} \label{eqn:ac-ohm-1} \\
        & S_{ji} = \left(Y_{ij} + Y^{c}_{ji} \right)^{*} \lvert V_{j} \rvert^{2} - Y^{*}_{ij} \frac{V_{i}^{*} V_{j}}{T^{*}_{ij}}, ~ \forall (i,j) \in \mathcal{E} \label{eqn:ac-ohm-2} \\
        & \sum_{\mathclap{k \in \mathcal{G}_{i}}} S^{g}_{k} - \sum_{\mathclap{k \in \mathcal{L}_{i}}} z^{d}_{k} S^{d}_{k} - \sum_{\mathclap{k \in \mathcal{H}_{i}}} z^{s}_{k} Y^{s}_{k} \lvert V_i \rvert^2 = \sum_{\mathclap{(i, j) \in \mathcal{E}_{i} \cup \mathcal{E}^{R}_{i}}} S_{ij}, ~ \forall i \in \mathcal{N} \label{eqn:ac-kcl} \\
        & \lvert S_{ij} \rvert \leq \overline{S}_{ij}, ~ S_{ij} \in \mathbb{C}, ~ \forall (i, j) \in \mathcal{E} \cup \mathcal{E}^R \label{eqn:power-magnitude-limit} \\
        & \underline{\theta}^{\Delta}_{ij} \leq \angle \left(V_{i} V^{*}_{j}\right) \leq \overline{\theta}^{\Delta}_{ij}, ~ \forall (i, j) \in \mathcal{E} \label{eqn:ac-phase-angle-difference-limits} \\
        & z^{v}_{i} \underline{V}_i \leq \lvert V_{i} \rvert \leq z^{v}_{i} \overline{V}_{i}, ~ V_{i} \in \mathbb{C}, ~ \forall i \in \mathcal{N} \label{eqn:ac-voltage-magnitude-bounds} \\
        & z^{g}_{i} \underline{S}^{g}_{i} \leq S^{g}_{i} \leq z^{g}_{i} \overline{S}^{g}_{i}, ~ S^{g}_{i} \in \mathbb{C}, ~ \forall i \in \mathcal{G} \label{eqn:ac-power-generation-bounds} \\
        & z^{v}_{i} \in \{0, 1\}, \; \forall i \in \mathcal{N}, \; z^{g}_{i} \in \{0, 1\}, \; \forall i \in \mathcal{G} \label{eqn:ac-discrete-indicators} \\
        & z^{d}_{i} \in [0, 1], \; \forall i \in \mathcal{L}, \; z^{s}_{i} \in [0, 1], \; \forall i \in \mathcal{H} \label{eqn:ac-continuous-indicators}
    \end{align}
    \label{eqn:power-constraints}
\end{subequations}
\caption{Power Network Modeling Requirements}
\label{model:ac-feasibility}
\end{model}

The MINLP formulation for AC power network feasibility, as defined for AC MLD analysis, is presented in Model \ref{model:ac-feasibility} and detailed by \cite{coffrin2018relaxations}.
Here, Constraints \eqref{eqn:ac-ohm-1} and \eqref{eqn:ac-ohm-2} model Ohm's law for lines, where $S_{ij} \in \mathbb{C}$ denotes the variable power along each line; $Y_{ij} \in \mathbb{C}$ and $Y_{ij}^{c} \in \mathbb{C}$ are constants denoting the line admittance and line charging; $V_{i} \in \mathbb{C}$ denotes the variable voltage at bus $i \in \mathcal{N}$; and $T_{ij} \in \mathbb{C}$ denotes constant transformer properties.
Constraints \eqref{eqn:ac-kcl} model power balances from Kirchhoff's current law for each bus, where $S_{k}^{g} \in \mathbb{C}$ denotes the variable power supplied by generator $k \in \mathcal{G}$; $S_{k}^{d} \in \mathbb{C}$ denotes the maximum power that can be delivered at load $k \in \mathcal{L}$; and $Y_{k}^{s}$ denotes the admittance of bus shunt $k \in \mathcal{H}$.
Note that $z_{k}^{d}$, $k \in \mathcal{L}_{i}$ allows each load to vary between zero and its predefined maximum, and $z_{k}^{s}$ allows for shedding fixed bus shunts from the network.
These modifications ensure that power balance constraints are satisfied in damaged networks.

Constraints \eqref{eqn:power-magnitude-limit}--\eqref{eqn:ac-continuous-indicators} impose engineering limits and variable bounds.
Constraints \eqref{eqn:power-magnitude-limit} bound the apparent power flow on each line, representing thermal limits.
Constraints \eqref{eqn:ac-phase-angle-difference-limits} ensure that each voltage phase angle difference is limited by predefined lower and upper bounds, $\underline{\theta}^{\Delta}_{ij}$ and $\overline{\theta}^{\Delta}_{ij}$, respectively.
Constraints \eqref{eqn:ac-voltage-magnitude-bounds} bound the voltage magnitude at each bus, where $\underline{V}_{i}$ and $\overline{V}_{i}$ denote lower and upper bounds, respectively.
Here, $z_{i}^{v} \in \{0, 1\}$ is a discrete variable that allows each bus to become de-energized when isolated from load or generation.
Similarly, Constraints \eqref{eqn:ac-power-generation-bounds} bound power generation, where $\underline{S}_{i}^{g}$ and $\overline{S}_{i}^{g}$ denote lower and upper bounds, respectively, and $z_{i}^{g} \in \{0, 1\}$ allows for each generator to become uncommitted when required to satisfy Constraints \eqref{eqn:ac-ohm-1} and \eqref{eqn:ac-ohm-2}.

\subsection{Natural Gas Transmission Network Modeling}
\label{subsection:natural-gas-modeling}

\paragraph{Notation for Sets}
A gas pipeline network is modeled using a directed graph $(\mathcal{J}, \mathcal{A})$, where $\mathcal{J}$ is the set of nodes (i.e., junctions) and $\mathcal{A}$ is the set of components that connect two nodes.
The sets of receipts (producers) and deliveries (consumers) are denoted by $\mathcal{R}$ and $\mathcal{D}$, respectively.
These components are considered to be attached to junctions $i \in \mathcal{J}$.
The subset of receipts attached to $i \in \mathcal{J}$ is denoted by $\mathcal{R}_{i}$ and the subset of deliveries by $\mathcal{D}_{i}$.
The sets of horizontal and short pipes are denoted by $\mathcal{P} \subset \mathcal{A}$ and $\mathcal{S} \subset \mathcal{A}$, respectively;
the set of resistors by $\mathcal{T} \subset \mathcal{A}$;
the set of valves and pressure-reducing regulators by $\mathcal{V} \subset \mathcal{A}$ and $\mathcal{W} \subset \mathcal{A}$, respectively;
and the set of compressors by $\mathcal{C} \subset \mathcal{A}$.
Additionally, the set of node-connecting components incident to $i \in \mathcal{J}$ where $i$ is the tail (respectively, head) of the arc is denoted by $\delta^{+}_{i} := \{(i, j) \in \mathcal{A}\}$ (respectively, $\delta^{-}_{i} := \{(j, i) \in \mathcal{A}\}$).
We next define the decision variables and constraints required to model a damaged gas network's steady-state operations.

\paragraph{Gas Network Modeling Requirements}
\begin{model}[t]
\caption{Gas Network Modeling Requirements}
\label{model:gas-feasibility}
\begin{subequations}
    \begin{align}
        & \sum_{\mathclap{(i, j) \in \delta^{+}_{i}}} f_{ij} - \sum_{\mathclap{(j, i) \in \delta^{-}_{i}}} f_{ji} = \sum_{\mathclap{k \in \mathcal{R}_{i}}} s_{k} - \sum_{\mathclap{k \in \mathcal{D}_{i}}} d_{k}, \; \forall i \in \mathcal{J} \label{eqn:mass-conservation} \\
        & p_{i}^{2} - p_{j}^{2} = w_{ij} f_{ij} \lvert f_{ij} \rvert, ~ \forall (i, j) \in \mathcal{P} \label{eqn:pipe-weymouth} \\
        & p_{i} - p_{j} = 0, ~ \forall (i, j) \in \mathcal{S} \label{eqn:short-pipe-pressure} \\
        & p_{i} - p_{j} = \tau_{ij} f_{ij} \lvert f_{ij} \rvert, ~ \forall (i, j) \in \mathcal{T} \label{eqn:resistor-darcy-weisbach} \\
        & \underline{f}_{ij} z_{ij} \leq f_{ij} \leq \overline{f}_{ij} z_{ij}, ~ z_{ij} \in \{0, 1\}, ~ \forall (i, j) \in \mathcal{V} \label{eqn:valve-flow-bounds} \\
        & p_{i} \leq p_{j} + (1 - z_{ij}) \overline{p}_{i}, ~ \forall (i, j) \in \mathcal{V} \label{eqn:valve-pressure-1} \\
        & p_{j} \leq p_{i} + (1 - z_{ij}) \overline{p}_{j}, ~ \forall (i, j) \in \mathcal{V} \label{eqn:valve-pressure-2} \\
        & \underline{f}_{ij} z_{ij} \leq f_{ij} \leq \overline{f}_{ij} z_{ij}, ~ z_{ij} \in \{0, 1\}, ~ \forall (i, j) \in \mathcal{W} \label{eqn:regulator-flow-bounds} \\
        & f_{ij} (p_{i} - p_{j}) \geq 0, ~ \forall (i, j) \in \mathcal{W} \label{eqn:regulator-direction} \\
        & \underline{\alpha}_{ij} p_{i} \leq p_{j} + (1 - z_{ij}) \underline{\alpha}_{ij} \overline{p}_{i}, ~ \forall (i, j) \in \mathcal{W} \label{eqn:regulator-pressure-1} \\
        & p_{j} \leq \overline{\alpha}_{ij} p_{i} + (1 - z_{ij}) \overline{p}_{j}, ~ \forall (i, j) \in \mathcal{W} \label{eqn:regulator-pressure-2} \\
        & \underline{\alpha}_{ij} p_{i} \leq p_{j} \leq \overline{\alpha}_{ij} p_{i}, ~ \forall (i, j) \in \mathcal{C} : \underline{f}_{ij} \geq 0 \label{eqn:compressor-pressures-1} \\
        & \underline{\alpha}_{ij} p_{i} \leq p_{j} \leq \overline{\alpha}_{ij} p_{i}, ~ \forall (i, j) \in \mathcal{C} : \underline{f}_{ij} < 0 \land \underline{\alpha}_{ij} = 1 \label{eqn:compressor-pressures-2-1} \\
        & f_{ij} (p_{i} - p_{j}) \leq 0, ~ \forall (i, j) \in \mathcal{C} : \underline{f}_{ij} < 0 \land \underline{\alpha}_{ij} = 1 \label{eqn:compressor-pressures-2-3} \\
        & y_{ij} \in \{0, 1\}, ~ \forall (i, j) \in \mathcal{C} : \underline{f}_{ij} < 0 \land \underline{\alpha}_{ij} \neq 1 \label{eqn:compressor-pressures-3-0} \\
        & p_{j} \leq \overline{\alpha}_{ij} p_{i} + (1 - y_{ij}) \overline{p}_{j}, ~ \forall (i, j) \in \mathcal{C} : \underline{f}_{ij} < 0 \land \underline{\alpha}_{ij} \neq 1 \label{eqn:compressor-pressures-3-1} \\
        & \underline{\alpha}_{ij} p_{i} \leq p_{j} + (1 - y_{ij}) \overline{p}_{i}, ~ \forall (i, j) \in \mathcal{C} : \underline{f}_{ij} < 0 \land \underline{\alpha}_{ij} \neq 1 \label{eqn:compressor-pressures-3-2} \\
        & p_{i} - p_{j} \leq y_{ij} \overline{p}_{i}, ~ \forall (i, j) \in \mathcal{C} : \underline{f}_{ij} < 0 \land \underline{\alpha}_{ij} \neq 1 \label{eqn:compressor-pressures-3-3} \\
        & p_{j} - p_{i} \leq y_{ij} \overline{p}_{j}, ~ \forall (i, j) \in \mathcal{C} : \underline{f}_{ij} < 0 \land \underline{\alpha}_{ij} \neq 1 \label{eqn:compressor-pressures-3-4} \\
        & \underline{f}_{ij} \leq f_{ij} \leq \overline{f}_{ij}, ~ \forall (i, j) \in \mathcal{A} \label{eqn:mass-flow-bounds} \\
        & 0 \leq \underline{p}_{i} \leq p_{i} \leq \overline{p}_{i}, ~ \forall i \in \mathcal{N} \label{eqn:pressure-bounds} \\
        & 0 \leq s_{k} \leq \overline{s}_{k}, ~ \forall k \in \mathcal{R}, ~ 0 \leq d_{k} \leq \overline{d}_{k}, ~ \forall k \in \mathcal{D} \label{eqn:sup-del-bounds}
    \end{align}
    \label{eqn:gas-constraints}
\end{subequations}
\end{model}

The MINLP formulation for gas network feasibility, as defined for MLD analysis, is presented in Model \ref{model:gas-feasibility} and detailed by \cite{tasseff2020natural}.
First, Constraints \eqref{eqn:mass-conservation} model nodal physics, i.e., mass flow conservation at junctions $i \in \mathcal{J}$.
Here, $f_{ij} \in \mathbb{R}$ denotes the variable mass flow along each node-connecting component; $s_{k} \in \mathbb{R}_{+}$ denotes the variable supply at receipt $k \in \mathcal{R}$; and $d_{k} \in \mathbb{R}_{+}$ denotes the variable demand (or load) at delivery $k \in \mathcal{D}$.

Constraints \eqref{eqn:pipe-weymouth}--\eqref{eqn:pressure-bounds} model the physics of node-connecting components.
Constraints \eqref{eqn:pipe-weymouth} model the Weymouth relationship for steady-state flow in a gas pipeline for each horizontal pipe $(i, j) \in \mathcal{P}$.
Here, $p_{i} \in \mathbb{R}_{+}$ denotes the variable pressure at junction $i \in \mathcal{J}$, and $w_{ij} \in \mathbb{R}_{+}$ denotes the constant mass flow resistance of the pipe.
These constraints are the most frequent sources of nonconvex nonlinearity in modeling the gas system.

Constraints \eqref{eqn:short-pipe-pressure} model short pipes in the network, which provide resistanceless mass transport between two junctions.
Constraints \eqref{eqn:resistor-darcy-weisbach} model resistors in the network, which act as surrogate components capable of modeling pressure losses elsewhere from pipes.
Here, pressure loss is modeled according to the Darcy-Weisbach equation, where $\tau_{ij} \in \mathbb{R}_{+}$ is the resistance, which is a function of the resistor's unitless drag factor and (possibly artificial) diameter.
Note that like Constraints \eqref{eqn:pipe-weymouth}, these constraints are also nonconvex nonlinear.

Constraints \eqref{eqn:valve-flow-bounds}--\eqref{eqn:valve-pressure-2} model valves in the network.
Here, the operating status of each valve $(i, j) \in \mathcal{V}$ is modeled using a discrete variable $z_{ij} \in \{0, 1\}$, where $z_{ij} = 1$ indicates an open valve and $z_{ij} = 0$ indicates a closed valve.
Constraints \eqref{eqn:valve-flow-bounds} prohibit flow across each valve when $z_{ij} = 0$.
Constraints \eqref{eqn:valve-pressure-1} and \eqref{eqn:valve-pressure-2} model, when a valve is open, that the pressures at connecting junctions are equal.
They also model the decoupling of junction pressures when the valve is closed.

Constraints \eqref{eqn:regulator-flow-bounds}--\eqref{eqn:regulator-pressure-2} model regulators (i.e., pressure-reducing valves) in the network.
Similar to valves, the status of each regulator is modeled using a discrete variable $z_{ij} \in \{0, 1\}$, where $z_{ij} = 1$ and $z_{ij} = 0$ indicate active and inactive statuses, respectively.
Constraints \eqref{eqn:regulator-flow-bounds} prohibit mass flow across each regulator when $z_{ij} = 0$.
Constraints \eqref{eqn:regulator-direction} ensure that mass flow across each regulator is in the same direction as the loss in pressure.
Constraints \eqref{eqn:regulator-pressure-1} and \eqref{eqn:regulator-pressure-2} model the remaining pressure dynamics.
Here, each regulator has a corresponding scaling factor, $\alpha_{ij}$.
This factor models the relationship between junction pressures when the regulator is active, i.e., $\alpha_{ij} p_{i} = p_{j}$.
The factor is further limited by the bounds $\underline{\alpha}_{ij} = 0 \leq \alpha_{ij} \leq \overline{\alpha}_{ij} = 1$.
Constraints \eqref{eqn:regulator-pressure-1} and \eqref{eqn:regulator-pressure-2} require that, when a regulator is active, pressures are defined according to the scaling relationship.
Otherwise, the pressures at the junctions connected by the regulator are decoupled.

Constraints \eqref{eqn:compressor-pressures-1}--\eqref{eqn:compressor-pressures-3-4} model compressors in the network.
Each compressor $(i, j) \in \mathcal{C}$ models an increase in pressure at junction $j \in \mathcal{J}$ by a variable scalar $\alpha_{ij}$.
Without loss of generality, bidirectional compression is not considered, although each compressor may allow for \emph{uncompressed} flow in the opposite direction.
These different behaviors of compressors are modeled by employing three different sets of constraints.
The first are Constraints \eqref{eqn:compressor-pressures-1} for compressors that \emph{prohibit} reverse flow, where $\underline{\alpha}_{ij}$ and $\overline{\alpha}_{ij}$ are minimum and maximum pressure ratios.
The second are Constraints \eqref{eqn:compressor-pressures-2-1} and \eqref{eqn:compressor-pressures-2-3} for compressors where reverse flow is \emph{allowed} and $\underline{\alpha}_{ij} = 1$.
Note that here, if $f_{ij} < 0$, then $p_{i} = p_{j}$.
Finally, Constraints \eqref{eqn:compressor-pressures-3-0}--\eqref{eqn:compressor-pressures-3-4} model compressors where uncompressed reverse flow is allowed and $\underline{\alpha}_{ij} \neq 1$.
In this case, the behavior of each compressor is disjunctive in its flow direction.
To model this disjunction, discrete variables $y_{ij} \in \{0, 1\}$ are introduced in Constraints \eqref{eqn:compressor-pressures-3-0} to model the direction of flow through each compressor.
Here, $y_{ij} = 1$ indicates flow from $i$ to $j$, and $y_{ij} = 0$ indicates flow from $j$ to $i$.
Constraints \eqref{eqn:compressor-pressures-3-1}--\eqref{eqn:compressor-pressures-3-4} model the pressures and pressure differences between junctions as per the specified flow direction and compression ratio bounds.

The remaining Constraints \eqref{eqn:mass-flow-bounds}--\eqref{eqn:sup-del-bounds} are variable bounds.
Constraints \eqref{eqn:mass-flow-bounds} are mass flow bounds, Constraints \eqref{eqn:pressure-bounds} are pressure bounds, and Constraints \eqref{eqn:sup-del-bounds} are receipt and delivery bounds.
Note that Constraints \eqref{eqn:sup-del-bounds} differ from the typical assumption of fixed supply and demand.
These modifications ensure mass conservation is satisfied in damaged networks.

\subsection{Interdependency Modeling}
\begin{figure}[t]
    \centering
    \includegraphics[width=1.00\linewidth]{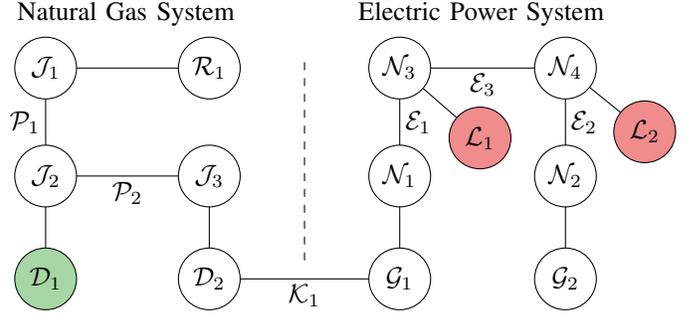}
    \caption{Diagrammatic representation of a small joint gas-power network. Here, $\mathcal{D}_{1}$ contributes to the objective term $\eta_{G}(\cdot)$ and $\mathcal{L}_{1}$, $\mathcal{L}_{2}$ contribute to $\eta_{P}(\cdot)$. Finally, the linkage between gas and power systems occurs at $\mathcal{K}_{1} = (\mathcal{D}_{2}, \mathcal{G}_{1})$.}
    \label{figure:joint-network}
\end{figure}

As in \cite{bent2018joint} and \cite{byeon2019unit}, gas and power systems are connected via heat rate curve models for gas-fired power generators, i.e.,
\begin{equation}
    \sum_{\mathclap{i : (i, j) \in \mathcal{K}}} h_{i}^{1} \Re(S_{i}^{g})^{2} + h_{i}^{2} \Re(S_{i}^{g}) + h_{i}^{3} z_{i}^{g} = d_{j}, ~ \forall j \in \mathcal{D}_{G} \label{eqn:heat-rate-constraints}.
\end{equation}
Each constraint links the real power generated at possibly multiple generators with a single gas delivery.
Here, $h_{i}^{1, 2, 3}$ are coefficients of the heat rate curve for $i \in \mathcal{G}$, and $\mathcal{K}$ is the set of linkages between gas-fired generators in $\mathcal{G}$ and their corresponding gas delivery points in $\mathcal{D}_{G} \subset \mathcal{D}$.
Furthermore, $h_{i}^{1} \geq 0$ for all $(i, j) \in \mathcal{K}$, and thus the left-hand side is always a convex function.
However, note that Constraint \eqref{eqn:heat-rate-constraints} is nonlinear nonconvex when $h_{i}^{1} \neq 0$.
Finally, the presence of $h_{i}^{3} z_{i}^{g}$ ensures that when $z_{i} = 0$, the intercept of the heat rate curve, and thus both generation and gas required, will be zero when a generator is uncommitted from the dispatch scenario.

A diagramatic illustration of the joint network model is illustrated in Figure \ref{figure:joint-network}.
Here, gas and power systems are linked by the single interdependency $\mathcal{K}_{1}$, which relates the delivery $\mathcal{D}_{2}$ to the generator $\mathcal{G}_{1}$.
Contributions to the gas and power delivery objectives, which are later described in Section \ref{subsection:mld-objective}, are represented by green and red colored nodes, respectively.

\subsection{Challenges}
Although independent gas and power MLD models were explored by \cite{tasseff2020natural} and \cite{coffrin2018relaxations}, respectively, the joint MLD problem that includes Constraints \eqref{eqn:power-constraints}{-}\eqref{eqn:heat-rate-constraints} is more challenging.
Most importantly, the nonlinear nonconvexities that appear in Models \ref{model:ac-feasibility} and \ref{model:gas-feasibility} arise primarily from different sources: Model \ref{model:ac-feasibility} includes many nonlinear equations with bilinear variable products, whereas Model \ref{model:gas-feasibility} includes more manageable quadratic nonlinear equations. 
To model them exactly, Model \ref{model:ac-feasibility} must be formulated as a highly challenging MINLP, but Model \ref{model:gas-feasibility} can be written as a more tractable mixed-integer nonconvex quadratic program.
These differences suggest potentially incompatible numerical methods and solving technologies in practice.

Although omitted here for brevity, this work also employed a number of important preprocessing steps used to ensure the construction of feasible damaged joint networks that satisfy Constraints \eqref{eqn:power-constraints}{-}\eqref{eqn:heat-rate-constraints}.
Finally, Models \ref{model:ac-feasibility} and \ref{model:gas-feasibility} constrain each system using steady-state physical assumptions.
In practice, modeling the transient dynamics of the gas system could be crucial.
However, as will be shown in subsequent sections, even the steady-state variant considered in this paper is computationally challenging.
This work is an important first step toward building MLD techniques that also consider transients.

\section{Maximal Load Delivery Formulations}
\label{section:mld-formulations}
This section derives the joint gas-power MLD formulations used throughout the remainder of this paper.
First, Section \ref{subsection:mld-objective} defines the competing objectives of the joint MLD problem.
Section \ref{subsection:mld-formulations} poses lexicographic and weighted MLD formulations that prioritize the gas-power delivery tradeoff in different ways.
Section \ref{subsection:micp-relaxation} derives MICP relaxations of the MINLP MLD formulations.
Finally, Section \ref{section:formulation-summary} summarizes the naming conventions used for these various MLD formulations, which are then empirically compared in Section \ref{section:computational-experiments}.

\subsection{Objectives of the Maximum Load Delivery Problem}
\label{subsection:mld-objective}
The objective of the MLD problem is to maximize the amount of \emph{nongeneration} gas and \emph{active} power load delivered simultaneously under a multi-contingency scenario.
Note that the maximization of \emph{nongeneration} gas, specifically, allows the model to decouple practical objectives of the gas system (e.g., delivery of fuel for residential heating) from practical objectives of the power system.
However, because the delivery of nongeneration gas load can inhibit the amount of active power generation, and thus active power delivered, there exists an important tradeoff between these two objectives.
For notational ease, we first write the two objective functions as
\begin{subequations}
\begin{align}
    \eta_{G}(d) &:= (\sum_{\mathclap{i \in \mathcal{D}^{\prime}}} \beta_{i} d_{i}) (\sum_{\mathclap{i \in \mathcal{D}^{\prime}}} \beta_{i} \overline{d}_{i})^{-1} \label{eqn:gas-objective} \\
    \eta_{P}(z^{d}) &:= (\sum_{\mathclap{i \in \mathcal{L}}} \beta_{i} z_{i}^{d} \lvert \Re({S}_{i}^{d}) \rvert) (\sum_{\mathclap{i \in \mathcal{L}}} \beta_{i} \lvert \Re({S}_{i}^{d}) \rvert)^{-1} \label{eqn:power-objective}.
\end{align}
\end{subequations}
Here, Equation \eqref{eqn:gas-objective} denotes the normalized sum of all prioritized nongeneration gas demand, where $\mathcal{D}^{\prime} := \mathcal{D} \setminus \{j : (i, j) \in \mathcal{K}\}$ (i.e., the set of all nongeneration gas deliveries), and $\beta_{i} \in \mathbb{R}_{+}$ is a predefined restoration priority for delivery $i \in \mathcal{D}^{\prime}$.
Similarly, Equation \eqref{eqn:power-objective} denotes the normalized sum of all prioritized active power loads.
Note that for all of the experiments considered in this study, $\beta_{i} = 1$ for all $i \in \mathcal{D}^{\prime} \cup \mathcal{L}$.

The tradeoff between nongeneration gas and active power load naturally lends the MLD problem to the broader category of multi-objective optimization.
A thorough survey of multi-objective optimization methods in engineering is presented by \cite{marler2004survey} and describes a number of techniques for specifying preferences among multiple objective functions.
These include weighted sum, weighted product, lexicographic, and bounded objective optimization methods.
In Section \ref{subsection:mld-formulations}, we define lexicographic and weighted sum variants of the MLD problem.

\subsection{Lexicographic and Weighted MLD Formulations}
\label{subsection:mld-formulations}
To explore the gas-power tradeoff, we introduce three MLD models that prioritize gas and power delivery in different ways.
The first is a lexicographic formulation that maximizes the amount of nongeneration gas load delivered \emph{first}.
This situation is representative of common contractual requirements for gas grid operators.
Here, the MLD is written as the program
\begin{equation}\tag{MLD-G}
\begin{aligned}
    & \text{maximize} & & \eta_{P}(z^{d}) \\
     & \text{subject to} & & \eta_{G}(d) \geq \eta_{G}(d^{*}) \\
     & & & \textnormal{Constraints} ~ \eqref{eqn:power-constraints}{-}\eqref{eqn:heat-rate-constraints},
\end{aligned}
\label{eqn:mld-gas}
\end{equation}
where $\eta_{G}(d^{*})$ is the optimal objective when maximizing gas delivery alone.
The second MLD is a similar formulation that maximizes the amount of active power load delivered first, i.e.,
\begin{equation}\tag{MLD-P}
\begin{aligned}
    & \text{maximize} & & \eta_{G}(d) \\
    & \text{subject to} & & \eta_{P}(z^{d}) \geq \eta_{P}(z^{d*}) \\
    & & & \textnormal{Constraints} ~ \eqref{eqn:power-constraints}{-}\eqref{eqn:heat-rate-constraints}.
\end{aligned}
\label{eqn:mld-power}
\end{equation}
The last is a single-level formulation that weights normalized sums of nongeneration gas and active power delivery, i.e.,
\begin{equation}\tag{MLD-W}
\begin{aligned}
    & \text{maximize} & & \lambda \eta_{G}(d) + (1 - \lambda) \eta_{P}(z^{d}) \\
    & \text{subject to} & & \textnormal{Constraints} ~ \eqref{eqn:power-constraints}{-}\eqref{eqn:heat-rate-constraints},
\end{aligned}
\label{eqn:mld-weighted}
\end{equation}
where $0 < \lambda < 1$ is a weighting parameter for the objective.

Note that \eqref{eqn:mld-gas}, \eqref{eqn:mld-power}, and \eqref{eqn:mld-weighted} are mixed-integer nonlinear, nonconvex programs.
The nonconvexities arise from three sources: (i) discrete operations of controllable components (e.g., $z_{i}^{g}$ for generator commitment); (ii) bilinear products that appear in both gas and power network physics (e.g., $V_{i} V_{j}^{*}$ in Ohm's law); and (iii) nonlinear equations used for satisfying physical relationships (e.g., the Weymouth equation for pipes).
In the following sections, we leverage a number of relaxations to render these problems more tractable.

\subsection{Relaxation of Bilinear Products and Nonlinear Equations}
\label{subsection:micp-relaxation}
\paragraph{Convexification of Power Physics}
The primary sources of nonconvexity in Model \ref{eqn:power-constraints} are the bilinear products that appear in Constraints \eqref{eqn:ac-ohm-1}--\eqref{eqn:ac-kcl} (e.g., $V_{i} V_{j}^{*}$).
A large body of literature has developed relaxations of similar terms, and for a comprehensive review, we refer the reader to \cite{molzahn2019survey}.
In this paper, we develop a model based on a second-order cone (SOC) relaxation of the AC power flow equations, first presented by \cite{jabr2006radial} and used for power MLD analysis in \cite{coffrin2018relaxations}.

The primary insight of the SOC formulation is that variable products ($\lvert V_{i} \rvert^{2}$ and $V_{i} V_{j}^{*}$) can be lifted into a higher-dimensional variable space ($W_{ii}$ and $W_{ij}$, respectively).
This renders terms involving these products linear, and the relaxation in the new $W$-space is ultimately strengthened via
\begin{equation}
    \lvert W_{ij} \rvert^{2} \leq W_{ii} W_{jj}, ~ \forall (i, j) \in \mathcal{E}.
\end{equation}
This is an SOC constraint, lending the formulation its name.

\paragraph{Convexification of Gas Physics}
Many nonconvexities in Model \ref{eqn:gas-constraints} appear in the form of nonlinear equations (e.g., Constraints \eqref{eqn:pipe-weymouth}) and bilinear variable products (e.g., Constraints \eqref{eqn:regulator-direction}).
To resolve both, direction variables $y_{ij} \in \{0, 1\}$ are first introduced for each node-connecting component $(i, j) \in \mathcal{A}$.
We also introduce variables $\pi_{i} \in \mathbb{R}_{+}$ to denote squared pressures $p_{i}^{2}$ for $i \in \mathcal{J}$.
This first allows for a partial linearization of the Weymouth equations for pipelines, i.e.,
\begin{equation}
    p_{i}^{2} - p_{j}^{2} = \pi_{i} - \pi_{j} = w_{ij} f_{ij} \lvert f_{ij} \rvert, ~ \forall (i, j) \in \mathcal{P} \label{eqn:pipe-weymouth-linearized}.
\end{equation}
Then, variables $\ell_{ij}$ for $(i, j) \in \mathcal{P}$ are introduced to model the difference in squared pressures across each pipe.
The introduction of $y$, $\pi$, and $\ell$, as well as convexly relaxing the equalities in Constraints \eqref{eqn:pipe-weymouth-linearized}, give rise to the convex relaxation
\begin{subequations}%
\begin{align}%
    & \pi_{j} - \pi_{i} \leq \ell_{ij} \leq \pi_{i} - \pi_{j}, ~ \forall (i, j) \in \mathcal{P} \label{eqn:weymouth-pipe-micqp-1} \\
    & \ell_{ij} \leq \pi_{j} - \pi_{i} + (2 y_{ij}) (\overline{\pi}_{i} - \underline{\pi}_{j}), ~ \forall (i, j) \in \mathcal{P} \label{eqn:weymouth-pipe-micqp-2} \\
    & \ell_{ij} \leq \pi_{i} - \pi_{j} + (2 y_{ij} - 2) (\underline{\pi}_{i} - \overline{\pi}_{j}), ~ \forall (i, j) \in \mathcal{P} \label{eqn:weymouth-pipe-micqp-3} \\
    & w_{ij} f_{ij}^{2} \leq \ell_{ij}, ~ \forall (i, j) \in \mathcal{P} \label{eqn:weymouth-pipe-micqp-4}.
\end{align}%
\label{eqn:weymouth-pipe-micqp}%
\end{subequations}
Note that Constraint \eqref{eqn:weymouth-pipe-micqp-4} is the primary physical relaxation, i.e., the Weymouth equation need not be satisfied with equality.

Convexification of the remaining nonlinear nonconvex terms in Model \ref{eqn:gas-constraints} is accomplished in a similar manner to the above.
Here, for brevity, we omit the derivation of the full mixed-integer convex relaxation used throughout the remainder of this study.
For a complete derivation and description of the relaxed mixed-integer convex model, we defer to \cite{tasseff2020natural}.

\paragraph{Convexification of Gas-fired Generation}
Constraints \eqref{eqn:heat-rate-constraints} are linear when $h_{i}^{1} = 0$ but nonconvex when $h_{i}^{1} > 0$.
In the latter case, Constraints \eqref{eqn:heat-rate-constraints} can be convexly relaxed as
\begin{equation}
    \sum_{\mathclap{i : (i, j) \in \mathcal{K}^{\prime}}} h_{i}^{1} \Re(S_{i}^{g})^{2} + h_{i}^{2} \Re(S_{i}^{g}) + h_{i}^{3} z_{i}^{g} \leq d_{j}, ~ \forall j \in \mathcal{D}_{G} \label{eqn:heat-rate-constraints-relaxed},
\end{equation}
where $\mathcal{K}^{\prime} := \{(i, j) \in \mathcal{K} : h_{i}^{1} \neq 0\}$.
However, in our experiments, all $h_{i}^{1}$ are zero, and the relaxation is not required.

\subsection{Summary of Formulations}
\label{section:formulation-summary}
The remainder of this paper compares two MLD formulations of Problems \eqref{eqn:mld-gas}, \eqref{eqn:mld-power}, and \eqref{eqn:mld-weighted}:
\begin{enumerate}
    \item (MLD-*): Exact MINLP formulations.
    \item (MLD-*-R): Formulations where power and gas constraints use SOC and MICP relaxations, respectively.
\end{enumerate}
These formulations provide different tradeoffs between model accuracy and computational performance.
An empirical evaluation of both allows us to quantify the effects of the relaxations, as well as to guide our subsequent MLD analyses.

\section{Computational Evaluation}
\label{section:computational-experiments}
In the following, Section \ref{subsection:experimental-setup} describes the networks, computational resources, and parameters used throughout the computational experiments;
Section \ref{subsection:n-k-experiments} compares the efficacy of exact and relaxed MLD formulations on randomized $N{-}k$ multi-contingency scenarios;
Section \ref{subsection:computational-performance} evaluates the runtime performance of formulations over the same experimental sets;
Section \ref{subsection:maximum_load_delivery_analysis} provides a proof-of-concept MLD analysis across the same experimental sets, illustrating the tradeoffs when lexicographically maximizing gas and power load delivery;
and Section \ref{subsection:pareto_analysis} provides a proof-of-concept Pareto analysis of load delivery on a single joint network.

\subsection{Benchmark Datasets and Experimental Setup}
\label{subsection:experimental-setup}

\begin{table*}[t]
    \begin{center}
        \caption{Summary of joint network datasets used in this study.}
        \label{table:networks}
        \begin{tabular}{|c|c!{\vrule width 1.5pt}c|c|c|c|c|c|c!{\vrule width 1.5pt}c|c|c|c|c!{\vrule width 1.5pt}c|}
            \hline
            Network & References & $\lvert \mathcal{J} \rvert$ & $\lvert \mathcal{P} \rvert$ & $\lvert \mathcal{S} \rvert$ & $\lvert \mathcal{T} \rvert$ & $\lvert \mathcal{V} \rvert$ & $\lvert \mathcal{W} \rvert$ & $\lvert \mathcal{C} \rvert$ & $\lvert \mathcal{N} \rvert$ & $\lvert \mathcal{E} \rvert$ & $\lvert \mathcal{G} \rvert$ & $\lvert \mathcal{L} \rvert$ & $\lvert \mathcal{H} \rvert$ & $\lvert \mathcal{K} \rvert$ \\ \hline
            \texttt{NG11-EP14}   & \cite{SABHJKKIOSSS17,babaeinejadsarookolaee2021power} & 11  & 8   & 0   & 0 & 1  & 0  & 2 & 14 & 20 & 5 & 11 & 1 & 1 \\ \hline
            \texttt{NG25-EP14}   & \cite{de2000gas,babaeinejadsarookolaee2021power,sanchez2016convex} & 25  & 24  & 0   & 0 & 0  & 0  & 6 & 14 & 20 & 5 & 11 & 1 & 2 \\ \hline
            \texttt{NG25-EP30}   & \cite{SABHJKKIOSSS17,babaeinejadsarookolaee2021power} & 25  & 19  & 1   & 1 & 0  & 2  & 3 & 30 & 41 & 6 & 21 & 2 & 1 \\ \hline
            \texttt{NG40-EP39}   & \cite{SABHJKKIOSSS17,babaeinejadsarookolaee2021power} & 40  & 39  & 0   & 0 & 0  & 0  & 6 & 39 & 46 & 10 & 21 & 0 & 4 \\ \hline
            \texttt{NG146-EP36}  & \cite{bent2018joint} & 146 & 93  & 0   & 0 & 0  & 42 & 29 & 36 & 121 & 91 & 35 & 2 & 34 \\ \hline
            \texttt{NG134-EP162} & \cite{SABHJKKIOSSS17,babaeinejadsarookolaee2021power} & 134 & 86  & 45  & 0 & 0  & 1  & 1 & 162 & 284 & 12 & 113 & 34 & 5 \\ \hline
            \texttt{NG135-EP179} & \cite{SABHJKKIOSSS17,babaeinejadsarookolaee2021power} & 135 & 141 & 0   & 0 & 0  & 0  & 29 & 179 & 263 & 29 & 104 & 40 & 12 \\ \hline
            \texttt{NG247-EP240} & \cite{price2011reduced} & 247 & 254 & 0   & 0 & 0  & 0  & 12 & 240 & 448 & 143 & 139 & 0 & 6 \\ \hline
            \texttt{NG603-EP588} & \cite{SABHJKKIOSSS17,babaeinejadsarookolaee2021power} & 603 & 278 & 269 & 8 & 26 & 44 & 5 & 588 & 686 & 167 & 379 & 68 & 12 \\ \hline
        \end{tabular}
    \end{center}
\end{table*}

The computational experiments in this paper consider gas and power networks of various sizes that appear in the literature or have been derived by subject matter experts.
These networks are summarized in Table \ref{table:networks}.
The networks in this table are named according to the number of junctions in the natural gas network (e.g., \texttt{NG11}) and the number of buses in the electric power network (e.g., \texttt{EP14}).
The references from which the gas, power, and/or joint network properties are derived appear in the second column of this table.
The numbers of natural gas and electric power system components of the joint networks vary substantially and are specified in the second and third delineated portions of Table \ref{table:networks}, respectively.
For networks that reference \cite{babaeinejadsarookolaee2021power}, heavily loaded variants of the corresponding electric power network datasets are used.

Joint gas-power network properties are summarized in the last column of Table \ref{table:networks}.
Here, \texttt{NG25-EP14} uses the linking and heat rate properties of the joint network instance developed by \cite{sanchez2016convex}, and \texttt{NG146-EP36} uses the properties of the instance developed by \cite{bent2018joint}.
Linkages within the \texttt{NG247-EP240} network were derived from open data, and heat rate curves were estimated in a manner similar to \cite{bent2018joint}.
The remaining networks combine instances from \textsc{GasLib} and \textsc{PGLib-OPF} to create new joint networks of various sizes.
The purpose of these new instances is twofold: (i) to explore the tractability of joint MLD instances as network sizes grow and (ii) to explore the tradeoffs involved in maximizing gas versus power delivery.
In these new instances, the number of gas-fired generators, $\lvert \mathcal{K} \rvert$, was estimated to be near $\min\{0.25 \lvert \mathcal{D} \rvert, 0.4 \lvert \mathcal{G} \rvert\}$, i.e., $\approx 25\%$ of all gas deliveries or $\approx 40\%$ of all generators.
After determining the total number of gas-fired generators, the largest-capacity generators in each power network were then assumed to be linked to the smallest-withdrawal delivery points in the gas network.
The heat rate at each gas-fired generator was then assumed to be equal to the proportion between the maximum withdrawal at the delivery point and the maximum power at the generator.
Note that these networks thus use synthetically generated linkages between \textsc{GasLib} and \textsc{PGLib-OPF} instances, and these linkages are not necessarily reflective of real-world datasets.
They are, however, instances where gas and power interdependencies are consequential, which in turn allows for a meaningful computational exploration of the MLD method.

All of the MLD formulations considered in this paper were implemented in the \textsc{Julia} programming language using the mathematical modeling layer \textsc{JuMP}, version 0.21 \cite{dunning2017jump}; version 0.9 of \textsc{GasModels}, a package for steady-state and transient natural gas network optimization \cite{gasmodels}; version 0.18 of \textsc{PowerModels}, a package for steady-state power network optimization \cite{coffrin2018powermodels}; and version 0.4 of \textsc{GasPowerModels}, a package for joint steady-state gas-power network optimization \cite{gaspowermodels}.
Furthermore, for the exact nonconvex nonlinear representation of Model \ref{model:ac-feasibility} in (MLD-*), the polar form of the AC power flow equations, introduced by \cite{carpentier1962contribution} and implemented by \cite{coffrin2018powermodels}, was used.
Similarly, for the exact representation of Model \ref{model:gas-feasibility}, the mixed-integer nonlinear nonconvex formulation described by \cite{tasseff2020natural} and implemented by \cite{gasmodels} was leveraged.

Each optimization experiment was prescribed a wall-clock limit of one hour on a node containing two Intel Xeon E5-2695 v4 processors, each with 18 cores @2.10 GHz, and 125 GB of memory.
For solutions of (MLD-W), version 0.7 of the open source \textsc{Juniper} MINLP solver was used \cite{kroger2018juniper}.
Within \textsc{Juniper}, \textsc{Ipopt} 3.12 was leveraged as the nonlinear programming solver, using a feasibility tolerance of $10^{-6}$ and the underlying linear system solver \textsc{MA57}, as recommended by \cite{tasseff2019exploring} for nonlinear network problems.
Note that \textsc{Juniper} does not provide global optimality guarantees for (MLD-W), and feasible solutions obtained from the solver serve only as \emph{lower bounds} on the true amount of maximum deliverable load.
For solutions of (MLD-*-R), \textsc{Gurobi} 9.1 was used with its default parameterization.
Here, since (MLD-*-R) is mixed-integer \emph{convex}, globally optimal solutions are obtained.
However, since (MLD-*-R) is a relaxation, a globally optimal solution corresponds only to an \emph{upper bound} on (MLD-*).

\subsection{Multi-contingency Damage Scenarios}
\label{subsection:n-k-experiments}
This section examines the robustness and accuracy of the exact and relaxed weighted MLD formulations, (MLD-W) and (MLD-W-R), respectively, with $\lambda = 0.5$.
Specifically, it studies these properties on large sets of randomized multi-contingency or $N{-}k$ scenarios, where $k$ indicates the number of components simultaneously removed from the joint gas-power network.
These scenarios are intended to capture the effects of severe multimodal network outages across joint systems.
In each scenario, a random selection of $15\%$ node-connecting components were assumed to be damaged (i.e., $k \approx 0.15 N$).
Through a parameter sensitivity study, we observed that this proportion of outages appeared to generate challenging MLD scenarios while providing interesting gas and power delivery tradeoffs among the coupled networks.
For each network, one thousand such scenarios were generated.

\begin{table}[t]
    \begin{center}
        \caption{Comparison of solver termination statuses over weighted objective MLD $N{-}k$ multi-contingency damage scenarios.}
        \label{table:n-k}
        \begin{tabular}{c|c|c|c|c|c|c|}
            \cline{2-7} & \multicolumn{3}{c|}{(MLD-W) $\%$ Cases} & \multicolumn{3}{c|}{(MLD-W-R) $\%$ Cases} \\
            \cline{1-7} \multicolumn{1}{|c|}{Network} & Conv. & Lim. & Inf. & Conv. & Lim. & Inf. \\
            \cline{1-7} \multicolumn{1}{|c|}{\texttt{NG11-EP14}} & $100.00$ & $0.00$ & $0.00$ & $100.00$ & $0.00$ & $0.00$ \\
            \cline{1-7} \multicolumn{1}{|c|}{\texttt{NG25-EP14}} & $99.90$ & $0.00$ & $0.10$ & $100.00$ & $0.00$ & $0.00$ \\
            \cline{1-7} \multicolumn{1}{|c|}{\texttt{NG25-EP30}} & $98.80$ & $0.10$ & $1.10$ & $99.70$ & $0.00$ & $0.30$ \\
            \cline{1-7} \multicolumn{1}{|c|}{\texttt{NG40-EP39}} & $99.00$ & $0.50$ & $0.50$ & $100.00$ & $0.00$ & $0.00$ \\
            \cline{1-7} \multicolumn{1}{|c|}{\texttt{NG146-EP36}} & $1.00$ & $83.70$ & $15.30$ & $100.00$ & $0.00$ & $0.00$ \\
            \cline{1-7} \multicolumn{1}{|c|}{\texttt{NG134-EP162}} & $26.70$ & $25.70$ & $47.60$ & $100.00$ & $0.00$ & $0.00$ \\
            \cline{1-7} \multicolumn{1}{|c|}{\texttt{NG135-EP179}} & $0.10$ & $95.20$ & $4.70$ & $100.00$ & $0.00$ & $0.00$ \\
            \cline{1-7} \multicolumn{1}{|c|}{\texttt{NG247-EP240}} & $0.00$ & $97.90$ & $2.10$ & $100.00$ & $0.00$ & $0.00$ \\
            \cline{1-7} \multicolumn{1}{|c|}{\texttt{NG603-EP588}} & $0.00$ & $70.30$ & $29.70$ & $100.00$ & $0.00$ & $0.00$ \\
            \cline{1-7}
        \end{tabular}
    \end{center}
\end{table}

Table \ref{table:n-k} compares statistics of solver termination statuses across all $N{-}k$ scenarios for each network and formulation.
Here, ``Conv.'' corresponds to the percentage of cases where the solver converged, ``Lim.'' to cases where the solver time or other solver limit was reached, and ``Inf.'' to cases that were classified as infeasible by the solver.
Although both formulations are typically capable of converging on cases containing tens of nodes, for larger networks, (MLD-W-R) clearly outperforms (MLD-W), solving nearly all $N{-}k$ instances.
The results are especially dramatic for the three largest networks, where only one of three thousand (MLD-W) cases converges but all (MLD-W-R) cases converge.
Note that three (MLD-W-R) cases are classified as infeasible due to numerical difficulties, but many more (MLD-W) cases are classified as infeasible due to the MINLP formulation and solver's greater tendency to converge to locally infeasible points.

\begin{table}[t]
    \begin{center}
        \caption{Comparison of solution quality for exact and relaxed joint MLD formulations over all $N{-}k$ contingency scenarios.}
        \label{table:exact-versus-relaxed}
        \begin{tabular}{c|c|c|c|c|}
            \cline{2-5} & \multicolumn{2}{c|}{(MLD-W) Solns.} & \multicolumn{2}{c|}{(MLD-W-R) Gap ($\%$)} \\
            \cline{1-5} \multicolumn{1}{|c|}{Network} & \# Compared & Mean Obj. & Mean & Median \\
            \cline{1-5} \multicolumn{1}{|c|}{\texttt{NG11-EP14}} & $1000$ & $0.61$ & $0.61$ & $0.03$ \\
            \cline{1-5} \multicolumn{1}{|c|}{\texttt{NG25-EP14}} & $999$ & $0.68$ & $1.45$ & $0.17$ \\
            \cline{1-5} \multicolumn{1}{|c|}{\texttt{NG25-EP30}} & $983$ & $0.45$ & $27.22$ & $0.09$ \\
            \cline{1-5} \multicolumn{1}{|c|}{\texttt{NG40-EP39}} & $990$ & $0.65$ & $0.45$ & $0.01$ \\
            \cline{1-5} \multicolumn{1}{|c|}{\texttt{NG146-EP36}} & $10$ & $0.75$ & $4.00$ & $1.01$ \\
            \cline{1-5} \multicolumn{1}{|c|}{\texttt{NG134-EP162}} & $267$ & $0.53$ & $52.63$ & $1.35$ \\
            \cline{1-5} \multicolumn{1}{|c|}{\texttt{NG135-EP179}} & $1$ & $0.59$ & $0.33$ & $0.33$ \\
            \cline{1-5} \multicolumn{1}{|c|}{\texttt{NG247-EP240}} & $0$ & -- & -- & -- \\
            \cline{1-5} \multicolumn{1}{|c|}{\texttt{NG603-EP588}} & $0$ & -- & -- & -- \\
            \cline{1-5}
        \end{tabular}
    \end{center}
\end{table}

Whereas Table \ref{table:n-k} measures the numerical reliability of exact and relaxed MLD formulations, Table \ref{table:exact-versus-relaxed} compares the solution quality of relaxed formulations with feasible lower bounds obtained from (MLD-W).
Here, ``$\#$ Compared'' corresponds to the number of cases used in each comparison, ``Mean Obj.'' is the mean objective value obtained by (MLD-W) over all compared instances, ``Mean'' is the mean relative gap between (MLD-W) and (MLD-W-R) objective values, and ``Median'' is the median relative gap between objective values.
In each such measurement, the relative gap is computed as
\begin{equation}
    \textnormal{Relative Gap} := \left(\frac{\tilde{\eta} - \eta}{\eta}\right) 100\%,
\end{equation}
where $\tilde{\eta}$ is the objective value obtained when solving (MLD-W-R) and $\eta$ is the objective value when solving (MLD-W).

We note that, for \texttt{NG25-EP30}, five instances were excluded in the comparison: the three infeasible (MLD-W-R) instances and two instances that implied a negative relative gap.
Proceeding with the analysis, the mean objective values for all sets of feasible solutions indicate that between around $50\%$ and $75\%$ of gas and power loads are being delivered across all multi-contingency scenarios.
Second, the mean relative gap between feasible solutions obtained by (MLD-W) and the upper bounds obtained by (MLD-W-R) are sometimes large, with the largest being $52.63\%$ across all \texttt{NG134-EP162} damage scenarios.

These extreme gaps have only two sources from which they can arise.
First, a feasible solution obtained by \textsc{Juniper} for an (MLD-W) instance is not guaranteed to be near the globally optimal solution.
That is, the globally optimal (MLD-W) objective value is potentially much larger than what \textsc{Juniper} reports at termination.
Second, since (MLD-W-R) is a relaxation, it \emph{upper-bounds} the globally optimal objective value of (MLD-W).
The median column in Table \ref{table:exact-versus-relaxed} reports measures of centrality without the outliers that are likely arising from the first source of discrepancy.
Through these measurements, (MLD-W-R) is observed to often provide reliable and tight bounds on the optimal objective of (MLD-W), with relative gaps most often ranging from nearly zero to less than $1.35\%$.
This indicates that the relaxation is capable of providing tight upper bounds on maximum capacities of damaged networks.

\subsection{Computational Performance}
\label{subsection:computational-performance}
\begin{figure}[t]

\centering
\subfloat{
    \centering
    \begin{tikzpicture}[baseline]
        \begin{axis}[legend cell align=left,enlargelimits=false,xtick pos=bottom,
                     ytick pos=left,height=3.75cm,ylabel=Cases solved,
                     scaled x ticks=false,xticklabel style={/pgf/number format/fixed},
                     width=0.4\linewidth,ymin=1,ymax=4100,restrict x to domain=0:0.250,xmin=0.0,xmax=0.250]
            \pgfplotstableread[col sep = comma]{data/performance-tens-FORMULATION.DWPGasModel-ACPPowerModel-build_mld_uc.csv}{\mldw};
            \pgfplotstableread[col sep = comma]{data/performance-tens-FORMULATION.CRDWPGasModel-SOCWRPowerModel-build_mld_uc.csv}{\mldwc};
            \addplot[very thick] table [x = solve_time, y = num_problems_solved]{\mldw};
            \addplot[very thick, dashed] table [x = solve_time, y = num_problems_solved]{\mldwc};
            \node[anchor=north west] at (rel axis cs:0.06,0.98) {\large{(S)}};
        \end{axis}
    \end{tikzpicture}
    \begin{tikzpicture}[baseline]
        \begin{semilogxaxis}[legend cell align=left,enlargelimits=false,
                             xtick pos=bottom,ytick pos=left,height=3.75cm,
                             ytick=\empty,width=0.7\linewidth,scaled ticks=false,
                             xlabel absolute,legend style={at={(0.595, 0.08)},anchor=south west,font=\tiny},
                             legend columns=1,ymin=1,ymax=4100,xmin=0.25,xmax=3600.0]
            \pgfplotstableread[col sep = comma]{data/performance-tens-FORMULATION.DWPGasModel-ACPPowerModel-build_mld_uc.csv}{\mldw};
            \pgfplotstableread[col sep = comma]{data/performance-tens-FORMULATION.CRDWPGasModel-SOCWRPowerModel-build_mld_uc.csv}{\mldwc};
            \addplot[very thick] table [x = solve_time, y = num_problems_solved]{\mldw};
            \addplot[very thick, dashed] table [x = solve_time, y = num_problems_solved]{\mldwc};
        \end{semilogxaxis}
    \end{tikzpicture}} \\
\subfloat{
    \centering
    \begin{tikzpicture}[baseline]
        \begin{axis}[legend cell align=left,enlargelimits=false,xtick pos=bottom,
                     ytick pos=left,height=3.75cm,ylabel=Cases solved,
                     scaled x ticks=false,xticklabel style={/pgf/number format/fixed},
                     width=0.4\linewidth,ymin=1,ymax=4100,restrict x to domain=0:0.25,xmin=0.0,xmax=0.25]
            \pgfplotstableread[col sep = comma]{data/performance-hundreds-FORMULATION.DWPGasModel-ACPPowerModel-build_mld_uc.csv}{\mldw};
            \pgfplotstableread[col sep = comma]{data/performance-hundreds-FORMULATION.CRDWPGasModel-SOCWRPowerModel-build_mld_uc.csv}{\mldwc};
            \addplot[very thick] table [x = solve_time, y = num_problems_solved]{\mldw};
            \addplot[very thick, dashed] table [x = solve_time, y = num_problems_solved]{\mldwc};
            \node[anchor=north west] at (rel axis cs:0.06,0.98) {\large{(M)}};
        \end{axis}
    \end{tikzpicture}
    \begin{tikzpicture}[baseline]
        \begin{semilogxaxis}[legend cell align=left,enlargelimits=false,
                             xtick pos=bottom,ytick pos=left,height=3.75cm,
                             ytick=\empty,width=0.7\linewidth,legend columns=1,
                             xlabel absolute,legend style={at={(0.96, 0.70)},anchor=north east,font=\footnotesize},
                             legend columns=1,ymin=1,ymax=4100,xmin=0.25,xmax=3600.0]
            \pgfplotstableread[col sep = comma]{data/performance-hundreds-FORMULATION.DWPGasModel-ACPPowerModel-build_mld_uc.csv}{\mldw};
            \pgfplotstableread[col sep = comma]{data/performance-hundreds-FORMULATION.CRDWPGasModel-SOCWRPowerModel-build_mld_uc.csv}{\mldwc};
            \addplot[very thick] table [x = solve_time, y = num_problems_solved]{\mldw};
            \addplot[very thick, dashed] table [x = solve_time, y = num_problems_solved]{\mldwc};
            \addlegendentry{(MLD-W)}
            \addlegendentry{(MLD-W-R)}
        \end{semilogxaxis}
    \end{tikzpicture}} \\
\subfloat{
    \centering
    \begin{tikzpicture}[baseline]
        \begin{axis}[legend cell align=left,enlargelimits=false,xtick pos=bottom,
                     ytick pos=left,height=3.75cm,ylabel=Cases solved,
                     scaled x ticks=false,xticklabel style={/pgf/number format/fixed},
                     width=0.4\linewidth,ymin=1,ymax=1025,restrict x to domain=0:0.25,xmin=0.0,xmax=0.25]
            \pgfplotstableread[col sep = comma]{data/performance-thousands-FORMULATION.DWPGasModel-ACPPowerModel-build_mld_uc.csv}{\mldw};
            \pgfplotstableread[col sep = comma]{data/performance-thousands-FORMULATION.CRDWPGasModel-SOCWRPowerModel-build_mld_uc.csv}{\mldwc};
            \addplot[very thick] table [x = solve_time, y = num_problems_solved]{\mldw};
            \addplot[very thick, dashed] table [x = solve_time, y = num_problems_solved]{\mldwc};
            \node[anchor=north west] at (rel axis cs:0.06,0.98) {\large{(L)}};
        \end{axis}
    \end{tikzpicture}
    \begin{tikzpicture}[baseline]
        \begin{semilogxaxis}[legend cell align=left,enlargelimits=false,
                             xtick pos=bottom,ytick pos=left,height=3.75cm,x tick label style={/pgf/number format/fixed},
                             xlabel={\hspace{-0.3\linewidth}Solution time (seconds)},
                             ytick=\empty,width=0.7\linewidth,
                             xlabel absolute,legend style={at={(0.97, 0.04)},anchor=south east},
                             legend columns=2,ymin=1,ymax=1025,xmin=0.25,xmax=3600.0]
            \pgfplotstableread[col sep = comma]{data/performance-thousands-FORMULATION.DWPGasModel-ACPPowerModel-build_mld_uc.csv}{\mldw};
            \pgfplotstableread[col sep = comma]{data/performance-thousands-FORMULATION.CRDWPGasModel-SOCWRPowerModel-build_mld_uc.csv}{\mldwc};
            \addplot[very thick] table [x = solve_time, y = num_problems_solved]{\mldw};
            \addplot[very thick, dashed] table [x = solve_time, y = num_problems_solved]{\mldwc};
        \end{semilogxaxis}
\end{tikzpicture}}
    \caption{Performance profiles comparing the efficiency of (MLD-W) and (MLD-W-R) over the $N{-}k$ instances described in Section \ref{subsection:n-k-experiments}. Here, the performance profiles are partitioned into three categories for (S) networks containing tens of nodes; (M) networks containing hundreds of nodes; and (L) networks containing more than a thousand nodes (i.e., \texttt{NG603-EP588}).}
    \label{figure:performance}
\end{figure}
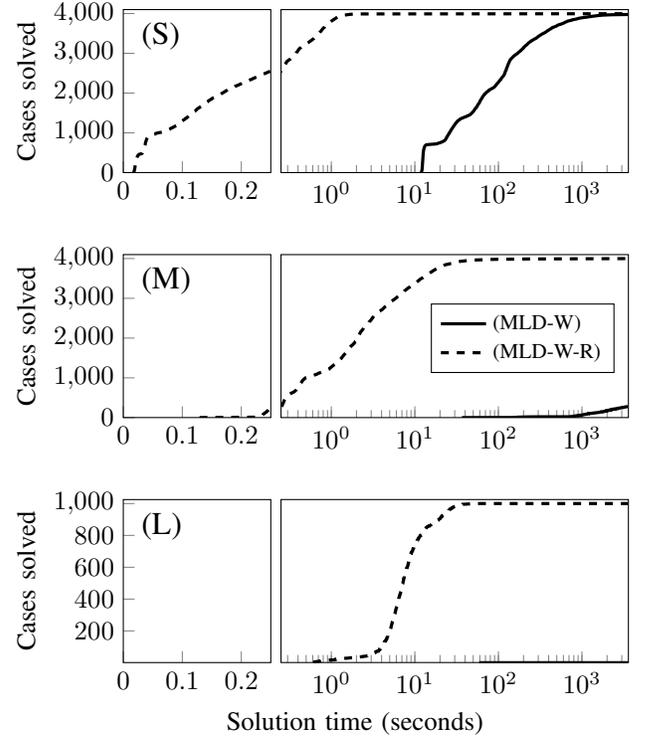

This section compares the performance of (MLD-W) and (MLD-W-R) using the instances described in Section \ref{subsection:n-k-experiments}.
The performance profiles for these cases are depicted in Figure \ref{figure:performance} and divided into three categories: (S) networks containing tens of nodes; (M) networks containing hundreds of nodes; and (L) networks containing more than a thousand nodes (i.e., \texttt{NG603-EP588}).
In all such categories, it is shown that the (MLD-W-R) formulation is able to solve substantially more problems than (MLD-W) in significantly shorter amounts of time.
For joint networks with tens of nodes, both formulations are able to solve many instances within the one hour time limit.
For networks with hundreds of nodes, (MLD-W-R) is capable of solving most instances within ten seconds, while (MLD-W) requires hundreds or thousands of seconds to solve only a small proportion.
For networks with thousands of nodes, (MLD-W-R) solves all instances within ten seconds, whereas (MLD-W) does not solve any.
The efficiency of (MLD-W-R) compared to (MLD-W) highlights its applicability to (i) real-time multi-contingency analysis and (ii) analyses that would require distributions of many multi-contingency scenarios.

\subsection{Proof-of-concept Maximum Load Delivery Analysis}
\label{subsection:maximum_load_delivery_analysis}
\begin{figure*}[!ht]
    \input{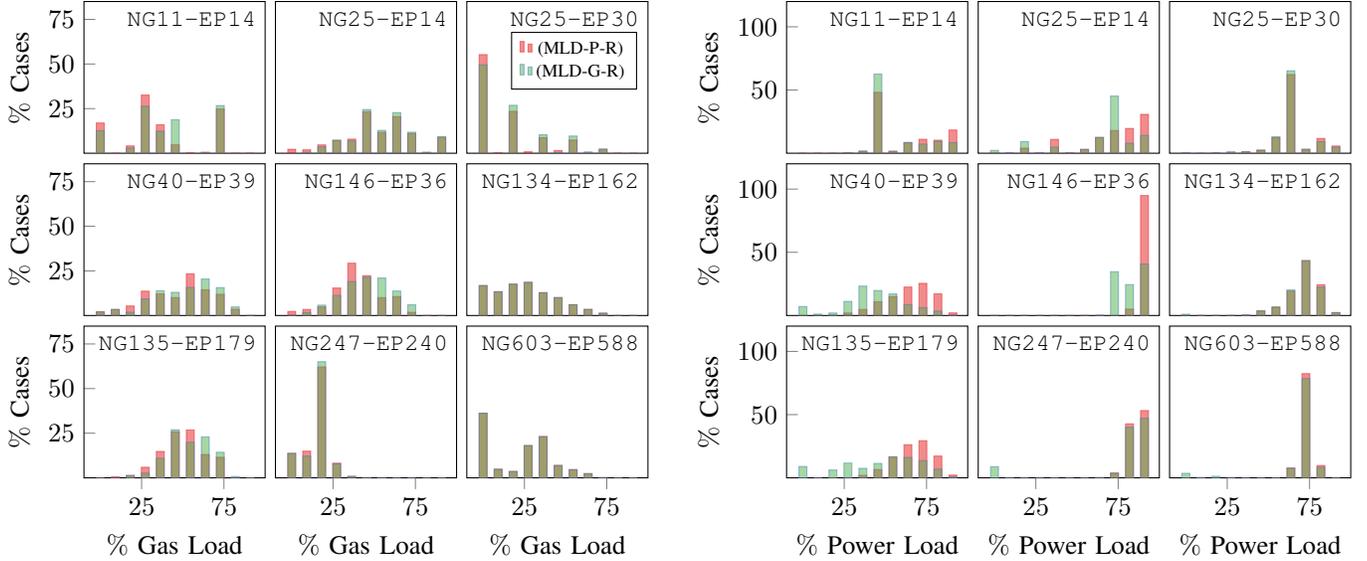}
    \caption{Histograms comparing the proportion of total gas and power load delivered across all solved $N{-}k$ scenarios for (MLD-G-R) and (MLD-P-R) variants. The $x$-axis indicates the proportion of load delivered, and the $y$-axis indicates the proportion of solved damage cases that deliver load within an interval.}
    \label{figure:n-k-histograms}
\end{figure*}

Whereas Sections \ref{subsection:n-k-experiments} and \ref{subsection:computational-performance} study the computational and accuracy tradeoffs between (MLD-W) and (MLD-W-R), this section provides a proof-of-concept MLD analysis using the (MLD-G-R) and (MLD-P-R) formulations on the same set of $N{-}k$ damage scenarios.
Figure \ref{figure:n-k-histograms} displays $18$ histograms that evaluate the proportions of gas and power loads delivered across \emph{solved} damage scenarios for the nine joint networks while using the two problem specifications.
Here, green bars correspond to histogram frequencies obtained from analyzing results of (MLD-G-R) solutions (i.e., gas prioritization) and red bars correspond to (MLD-P-R) solutions (i.e., power prioritization).
Brown, overlapping bars correspond to frequencies that appear in both (MLD-P-R) and (MLD-G-R) histograms.
These results indicate qualitative differences in the hypothetical robustness of each joint network.
They also display the extremal tradeoffs between prioritizing gas versus power delivery in the presence of extreme outages.
Finally, they indicate the sensitivity of each gas or power network to the interdependencies that link them.
These histograms serve as basic proofs of concept for real-world MLD analyses.

The left half of Figure \ref{figure:n-k-histograms} displays histograms of maximum gas load delivered in the presence of severe $N{-}k$ outages.
First, note that these histograms display a variety of load distributions across the cases and networks considered.
Some networks, e.g., \texttt{NG25-EP30}, \texttt{NG247-EP240}, and \texttt{NG603-EP588} suggest gas grids that are highly sensitive to the outages considered, with large proportions of damaged networks often incapable of delivering more than $50\%$ of gas load.
Other networks, e.g., \texttt{NG40-EP39}, \texttt{NG146-EP36}, and \texttt{NG135-EP179} show less severe but still substantial sensitivities to these outages.
The remaining networks display gas network sensitivities somewhere between these two extremes.

The overlapping histograms also display the tradeoffs encountered when prioritizing gas versus power delivery.
In the three joint networks \texttt{NG25-EP14}, \texttt{NG134-EP162}, and \texttt{NG603-EP588}, gas and power interdependencies are mostly inconsequential, and prioritizing either gas or power barely affects the maximum gas capacity.
This is likely a result of excess generation capacity in the corresponding power networks.
Other networks, e.g., \texttt{NG40-EP39}, \texttt{NG146-EP36}, and \texttt{NG135-EP179} show more interesting tradeoffs, where prioritizing either gas or power results in substantial changes in the overall maximum load distributions.
The remaining networks show less interesting tradeoffs, although \texttt{NG11-EP14} displays large tradeoffs, likely due to the drastic effects that even minor outages can have on the relatively small network.

The right half of Figure \ref{figure:n-k-histograms} displays histograms of maximum active power delivered in the presence of the $N{-}k$ outages.
First, the four networks \texttt{NG25-EP30}, \texttt{NG134-EP162}, \texttt{NG247-EP240}, and \texttt{NG603-EP588} appear robust to outages in the joint network and are often capable of delivering more than $75\%$ of the original power load.
The remaining networks see a greater variety in their maximum load distributions.
Whereas some networks, e.g., \texttt{NG25-EP30}, \texttt{NG134-EP162}, and \texttt{NG603-EP588}, appear less reliant on gas-fired power generators, the remaining networks exhibit more drastic changes when prioritizing gas versus power delivery.
The most extreme example appears to be \texttt{NG146-EP36}, which is often capable of delivering a large amount of power across all $N{-}k$ cases when power is prioritized but also often loses more than $25\%$ capacity when gas delivery is prioritized.

We remark that, to solve (MLD-G-R) and (MLD-P-R), inner- and outer-level problems of the lexicographic maximization are solved sequentially.
For example, to solve (MLD-G-R), (i) the inner level problem maximizing $\eta_{G}(d)$ is solved, yielding a solution $d^{*}$, then (ii) $\eta_{P}(z^{d})$ is maximized, subject to Constraints \eqref{eqn:power-constraints}{-}\eqref{eqn:heat-rate-constraints} and $\eta_{G}(d) \geq \eta_{G}(d^{*}) - \epsilon$.
The latter ensures that nongeneration gas load delivered in the outer-level is at least that of the inner level, minus some feasibility tolerance $\epsilon$, taken in this study to be $10^{-7}$.
A similar algorithm is used for (MLD-P-R).
We note that the general algorithm is not as numerically reliable as (MLD-W-R) and does not solve $469$ of the $18{,}000$ $N{-}k$ cases considered in this subsection.
This could be alleviated with a larger $\epsilon$ or direct use of lexicographic features available in some solvers (e.g., \textsc{Gurobi}).

\subsection{Proof-of-concept Pareto Analysis}
\label{subsection:pareto_analysis}
Together, (MLD-G-R), (MLD-P-R), and (MLD-W-R) allow for a variety of prioritizations of gas versus power load.
As such, they serve as powerful tools for exploring the wide range of possible MLD solutions based on the relative importance of gas versus power delivery.
This can provide gas and power grid managers with best-case capacity estimates depending on the type of coordination between the two systems.
In turn, this enables a better understanding of the extremely complex yet practically important tradeoffs encountered during the operation of a damaged joint network.
Whereas Sections \ref{subsection:n-k-experiments} through \ref{subsection:maximum_load_delivery_analysis} focus on analyzing performance and qualitative aspects of MLD analyses across a large number of joint networks, this section focuses on providing a proof-of-concept Pareto analysis on a single joint network, \texttt{NG146-EP36}.

\begin{figure}
    \centering
\begin{tikzpicture}
    \begin{axis}[height=5.00cm,enlargelimits=false,width=\linewidth,
                 ylabel=Active power load ($\%$),
                 xlabel=Nongeneration gas load ($\%$),
                 every node near coord/.append style={xshift=-1.5em,yshift=-2.0em},
                 legend style={at={(0.97, 0.80)},font=\tiny},anchor=north east,
                 cycle list name=Set1-9,xtick pos=bottom,ytick pos=left]
        \addplot [mark=*, dashed] table [x=percent_gas_served, y=percent_active_power_served, col sep=comma] {data/pareto.csv};
        \addplot [nodes near coords, point meta=explicit symbolic, only marks] table [x=percent_gas_served, y=percent_active_power_served, meta=weight, col sep=comma] {data/pareto.csv};
    \end{axis}
\end{tikzpicture}
    \caption{Pareto front approximation of total active power load versus nongeneration gas load delivered over one thousand \texttt{NG146-EP36} $N{-}k$ scenarios.}
    \label{figure:pareto}
\end{figure}
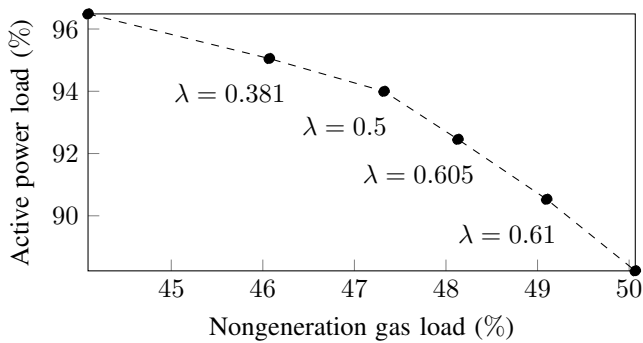

Figure \ref{figure:pareto} shows a linearly-interpolated approximation of the Pareto front for mean active power versus gas delivery across the same set of $N{-}k$ scenarios considered in previous sections.
Here, the upper-left and lower-right endpoints correspond to means obtained from the (MLD-P-R) and (MLD-G-R) problem formulations, respectively.
Interior data points correspond to means obtained from the (MLD-W-R) formulation, where the tradeoff parameter $\lambda$ was varied to determine interesting and distinct points on the Pareto front.

First, note that when prioritizing power delivery, on average, $96\%$ of active power is delivered but less than $45\%$ of nongeneration gas is delivered.
When gas delivery is prioritized, $88\%$ of power is delivered, while $50\%$ of gas is delivered.
Between these two extremes, the amount of gas and power increases and decreases, respectively, with increases in $\lambda$.
For $\lambda \lessapprox 0.5$, active power decreases more slowly as a function of $\lambda$, and for $\lambda \gtrapprox 0.5$, the rate of decrease appears larger.
In this case, $\lambda \approx 0.5$ happens to represent a value where (MLD-W-R) begins to prefer maximization of gas delivery over power delivery.
Thus, in practice, a point near this value of $\lambda$ could be one which maximizes simultaneous delivery of the two quantities while having practically equal prioritizations.

\section{Conclusion}
\label{section:conclusion}
Recent increases in gas-fired power generation have amplified interdependencies between natural gas and power transmission systems.
These interdependencies have engendered greater vulnerabilities to gas and power grids, where natural or man-made disruptions can require the curtailment of load in one or both systems.
To address the challenge of estimating maximum joint network capacities under these disruptions, this study considered the task of determining feasible steady-state operating points for severely damaged joint networks while ensuring the maximal delivery of gas and power loads simultaneously.
Mathematically, this task was represented as the mixed-integer, nonlinear nonconvex joint MLD problem.

Three variants of the MLD problem were formulated: one that prioritizes gas delivery, one that prioritizes power delivery, and one that assumes a linear tradeoff between the two objectives.
To increase the tractability of these problems, a mixed-integer convex relaxation of the joint network's physical constraints was proposed.
To demonstrate the relaxation's effectiveness, exact and relaxed MLD formulations were computationally compared across a variety of $N{-}k$ scenarios.
The relaxation was found to be a fast and reliable means for determining bounds on capacities of damaged networks.

Two proofs of concept were then provided to showcase the analytical power of the relaxed MLD problems.
The first provided comparisons between prioritizing gas versus power delivery in an MLD analysis.
These examples showcased the sometimes substantial tradeoffs that should be considered in extreme outage scenarios.
The second proof of concept provided a Pareto front approximation of gas versus power delivery across $N{-}k$ scenarios using a single joint network.
These proofs of concept highlight that the efficacy of the relaxation-based MLD method makes it a potentially valuable tool for complex real-world decision support applications.

Future work will focus on extending the MLD approaches developed in this paper.
First, additional gas and power relaxations should be considered to more accurately and efficiently scale to joint networks containing many thousands of nodes.
Preprocessing routines, such as optimization-based bound tightening, may also aid in improving existing relaxations.
Second, the current problem assumes the full coordination between gas and power systems when deciding operations that maximize load delivery.
The modeling of bidding mechanisms that drive both systems could provide more accurate joint capacity estimates.
Finally, capturing transient dynamics in gas networks is sometimes crucial for understanding the effects of network disruptions, which may only be realized long after the disruption occurs.
Future work should consider these transient effects when modeling load delivery in the gas network.

\section*{Acknowledgments}
The authors gratefully acknowledge Drs. David Fobes and Kaarthik Sundar for their contributions to the \textsc{InfrastructureModels} software packages.
They also thank the administrators of the Darwin computing cluster at Los Alamos National Laboratory.
This work was conducted under the auspices of the National Nuclear Security Administration of the U.S. Department of Energy at Los Alamos National Laboratory under Contract No. 89233218CNA000001.

\bibliographystyle{IEEEtran}
\bibliography{bibliography}

\vspace{-3em}
\begin{IEEEbiographynophoto}{Byron Tasseff}
received the B.S. degree in physics from the University of Northern Iowa, Cedar Falls, IA, USA, in 2012 and the M.S. and Ph.D. degrees in industrial and operations engineering from the University of Michigan, Ann Arbor, MI, USA, in 2018 and 2021, respectively.
He is currently a staff scientist at Los Alamos National Laboratory, where his research interests involve developing optimization techniques for problems involving fluids and critical infrastructure (e.g., water and gas networks).
\end{IEEEbiographynophoto}
\vspace{-3em}
\begin{IEEEbiographynophoto}{Carleton Coffrin}
received the B.S. degree in computer science and the B.F.A. degree in theatrical design from the University of Connecticut, Storrs, CT, USA, in 2006 and the M.S. and Ph.D. degrees from the Brown University, Providence, RI, USA, in 2010 and 2012, respectively.
He is currently a staff scientist at Los Alamos National Laboratory, where he studies the application of optimization methods to problems involving infrastructure networks.
\end{IEEEbiographynophoto}
\vspace{-3em}
\begin{IEEEbiographynophoto}{Russell Bent}
received the Ph.D. degree in computer science from the Brown University, Providence, RI, USA, in 2005.
He is currently a staff scientist at Los Alamos National Laboratory Applied Mathematics and Plasma Physics group.
He is also an Associate Editor for the INFORMS Journal of Computing.
At Los Alamos National Laboratory, he leads a team of researchers focused on developing next-generation algorithms for planning, operating, and designing critical infrastructure systems.
\end{IEEEbiographynophoto}

\end{document}